\documentclass[10pt]{article}
\usepackage{amsfonts}
\usepackage{amsmath}
\usepackage{mathrsfs,url}
\usepackage{mathrsfs, amscd,amssymb,amsthm,amsmath,bm,graphicx,psfrag,subfigure}

\makeatletter

\renewcommand{\@seccntformat}[1]{{\csname the#1\endcsname}{\normalsize .}\hspace{.5em}}
\makeatother

\def \[{\begin{equation}}
\def \]{\end{equation}}
\def \per{{\rm per}}
\newtheorem{thm}{Theorem}[section]

\newtheorem{defi}[thm]{Definition}
\newtheorem{claim}{Claim}

\newtheorem{lem}[thm]{Lemma}

\newenvironment{wst}
{\setlength{\leftmargini}{1.5\parindent}
 \begin{itemize}
 \setlength{\itemsep}{-1.1mm}}
{\end{itemize}}

\begin{document}
\setlength{\baselineskip}{15pt}
\begin{center}{\Large \bf
Edge-grafting theorems on permanents of the Laplacian matrices of graphs and their applications\footnote{Financially supported by the National Natural Science
Foundation of China (Grant No. 11071096)  and the Special Fund for Basic Scientific Research of Central Colleges (CCNU11A02015).}}

\vspace{2mm}

{\large Shuchao Li\footnote{E-mail: lscmath@mail.ccnu.edu.cn (S.C.
Li)}, Yan Li}\vspace{2mm}

Faculty of Mathematics and Statistics,  Central China Normal
University, Wuhan 430079, P.R. China\vspace{1mm}

\end{center}

\noindent {\bf Abstract}:\ The trees, respectively unicyclic graphs, on
$n$ vertices with the smallest Laplacian permanent are studied. 
In this paper, by edge-grafting transformations,
the $n$-vertex trees of given bipartition having the second and third smallest Laplacian permanent are identified.
Similarly, the $n$-vertex bipartite unicyclic graphs of given bipartition having the
first, second and third smallest Laplacian permanent are characterized. Consequently, the $n$-vertex bipartite unicyclic
graphs with the first, second and third smallest Laplacian permanent are determined.

\vspace{2mm} \noindent{\it Keywords}: Laplacian matrix; Laplacian coefficient; Permanent; Tree; Unicyclic graph; Bipartition

\vspace{2mm}


{\setcounter{section}{0}

\section{\normalsize Introduction}
Let $G=(V_G,E_G)$ be a simple connected graph with vertex set
$V_G=\{v_1,\ldots,v_n\}$ and edge set $E_G\not=\emptyset$. The
adjacency matrix $A(G)=(a_{ij})$ of $G$ is an $n\times n$ symmetric
matrix with $a_{ij}=1$ if and only if $v_i, v_j$ are adjacent and 0
otherwise. Since $G$ has no loops, the main diagonal of $A(G)$
contains only 0's. Denote the degree of $v_i$ by $d_G(v_i)$ (or
$d_i$) for $i=1,\ldots,n$, and let $D(G)$ be the diagonal matrix
whose $(i,i)$-entry is $d_i, i=1,2,\ldots, n$. The matrix
$L(G)=D(G)-A(G)$ is called the \textit{Laplacian matrix} of $G$. Of
course, $L(G)$ depends on the ordering of the vertices of $G$.
However, a different ordering leads
to a matrix which is permutation similar to $L(G)$. 
The matrix $Q(G)=D(G)+A(G)$ has been called the
\textit{signless Laplacian matrix} of $G$. 
For survey papers on this matrix the reader is referred to \cite{C-R-S2,C-R-S1,C-R-S3}.

If the vertex set of the connected graph $G$ on $n$ vertices can be partitioned into two subsets $V_1$ and $V_2$ such
that each edge joins a vertex of $V_1$ to a vertex of $V_2$, then $G$ has a $(p,
q)$-\textit{bipartition} where $|V_1|=p$ and $|V_2|=q.$ Without loss of generality we may assume that $p \leq q.$

A connected graph with $n$
vertices and $n$ edges is called a \textit{unicyclic graph}. For
convenience, let $\mathscr{T}_n^{p,q}$ (resp. $\mathscr{U}_n^{p,q}$)
be the set of all $n$-vertex trees (resp.
bipartite unicyclic graphs) with a $(p, q)$-bipartition, and let $\mathscr{U}_n$ be the
set of all bipartite unicyclic graphs on $n$ vertices.

Throughout we denote by $P_n,\, S_n$ and $C_n$ the path, star and cycle on
$n$ vertices, respectively. $G-v, G-uv$ denote the graph obtained from $G$ by deleting vertex $v \in V_G$,
or edge $uv \in E_G$, respectively (this notation is naturally extended if more than one vertex or edge is
deleted). Similarly, $G + uv$ is obtained from $G$ by adding edge $uv \not\in E_G$. The \textit{distance} between vertices
$u$ and $v$ in $G$ is denoted by $d_G(u,v)$.
Let $PV(G)$ denote the set of all pendant vertices of $G.$


The \textit{permanent} of $X = (x_{ij})\in M_{n\times n}$, denoted by $\per X$, is the quantity
$$\per X = \sum_{\sigma\in S_n}\prod_{t=1}^nx_{t\sigma(t)},$$
where $S_n$ is the symmetric group of degree $n$; see \cite{H-M}. It was suggested in
\cite{B-Z6} to use the polynomial \per$(xI-L(G))$ to distinguish
non-isomorphic trees. For more progress on the quantity \per$(\cdot)$,
the reader may be referred to \cite{F-Z-Z}.

The first research paper on permanent of the Laplacian matrix was
\cite{B-Z6}, in which lower bounds for the permanent of $L(G)$ were
conjectured by Merris, Rebman and Watkins. These lower bounds on
\per$L(G)$ were proved by Brualdi and Goldwasser \cite{B-Z1} and
Merris \cite{B-Z06}. For more recent results on Laplacian (resp. signless
Laplacian) permanent one is referred to  \cite{G-H-L,Z-L1,Z-L2}.

The Laplacian polynomial $\mu(G, \lambda)$ of $G$ is the \textit{characteristic polynomial} of its Laplacian matrix
$L(G)$, that is,
$$
\mu(G, \lambda) = \det(\lambda I_n - L(G)) = \sum_{k=0}^n(-1)^kc_k\lambda^{n-k}.
$$
It is easy to see that $c_0(G) = 1, c_1(G) = 2|E_G|, c_n(G) = 0, c_{n-1}(G) = n\tau(G)$, where $\tau(G)$ denotes
the number of spanning trees of $G$. For two $n$-vertex graphs $G_1$ and $G_2$, we say that $G_1$ is dominated by $G_2$ and write $G_1 \preceq G_2$, if
$c_k(G_1) \leq c_k(G_2)$ holds for all Laplacian coefficients $c_k, k = 0, 1, \ldots, n$. If $G_1 \preceq G_2$ and there exists $j$ such that $c_j(G_1)<c_j(G_2)$, then we write $G_1 \prec G_2.$

Note that the Laplacian coefficients have combinatorial significant, hence the research on the
Laplacian coefficients of graphs has received great attention in recent years; see \cite{S-S-H,6,7,8,9,12,14,15}
and the references therein. Zhou and Gutman \cite{B-Z-G} showed that among all trees of order $n$, the $k$th coefficient $c_k$ is the largest
when the tree is a path and is the smallest for a star, $k=0, 1, \ldots, n$. In view of Theorems 2.4 and 2.5 in \cite{B-Z1} the counterparts of these results for
the Laplacian permanent of trees are as the following.
\begin{thm}[\cite{B-Z1}]
Let $T$ be a tree with $n$ vertices. Then
$$
2(n-1)\le \per\ L(T)\le \frac{2-\sqrt{2}}{2}(1+\sqrt{2})^n+\frac{2+\sqrt{2}}{2}(1-\sqrt{2})^n.
$$
The left equality holds if and only if $T$ is a star, whereas the right equality holds if and only if
$T$ is a path.
\end{thm}

Brualdi and Goldwasser \cite{B-Z1} showed that $T_{n,m}$ is the unique tree among the $n$-vertex trees each of which contains an $m$-matching
having the minimum Laplacian permanent, where $T_{n,m}$ is the tree obtained from the star graph $S_{n-m+1}$ by attaching a pendant edge to each of certain $m-1$ non-central vertices of $S_{n-m+1}.$ Ili\'c \cite{7} showed that $T_{n,m}$ is also the unique $n$-vertex tree with
given matching number $m$ which simultaneously minimizes all the Laplacian coefficients.
It is then natural to conjecture that among the class of graphs, a particular graph has the smallest Laplacian permanent, then that particular graph also minimizes  all of its Laplacian coefficients in that class of graphs, and vice versa. This mathematical phenomenon
is further studied in \cite{B-Z1, 6}. We know from \cite{6} that among the $n$-vertex trees of diameter $d$, caterpillar $T_{n,d, \lfloor d/2\rfloor}$ has the minimum Laplacian coefficient $c_k$, for every $k = 0, 1, \ldots, n,$ whereas we know from \cite{B-Z1} that
among the $n$-vertex trees of diameter $d$, the broom $T_{n,d,2}$ has the minimum Laplacian permanent. Graphs
$T_{n,d, \lfloor d/2\rfloor}$ and $T_{n,d,2}$ are depicted in Fig. 1. This implies that there is no monotone relationship
between the Laplacian coefficients and the Laplacian permanent of graphs. Yet we lack a better
understanding of this relationship.
\begin{figure}[h!]
\begin{center}
\psfrag{1}{$v_1$}\psfrag{2}{$v_2$} \psfrag{3}{$v_3$}\psfrag{5}{$v_{d}$}\psfrag{6}{$v_{d+1}$}\psfrag{4}{$v_{\lfloor d/2\rfloor}$}
\psfrag{7}{$u_1$}\psfrag{8}{$u_2$} \psfrag{9}{$u_{n-d-1}$}\psfrag{a}{$u_{n-d-1}$}\psfrag{A}{$T_{n,d, \lfloor d/2\rfloor}$}
\psfrag{B}{$T_{n,d,2}$} \psfrag{g}{$G$}\psfrag{h}{$G[u\rightarrow w;2]$}
\psfrag{f}{$u_2$} \psfrag{i}{$u_2$}
\includegraphics[width=110mm]{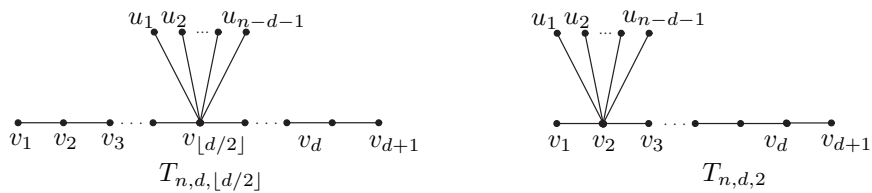}\\
\caption{Graphs $T_{n,d, \lfloor d/2\rfloor}$ and $T_{n,d,2}$.}
\end{center}
\end{figure}

An interesting fact is that among $n$-vertex trees with a given
bipartition, the extremal one that minimizes the Laplacian permanent \cite{B-Z1} is exactly
the one that simultaneously minimizes all Laplacian coefficients; see \cite{W-G-Y}.
Up to now it is natural for us to find further examples, where, like in Theorem 1.1, if the Laplacian permanent is minimized (maximized) by a
particular graph in a class of graphs, then that particular graph minimizes (maximizes) all the Laplacian coefficients in that class of
graphs.

Motivated by \cite{B-Z1, W-G-Y}, in this paper, we use a new and unified method to show some known results on the Laplacian permanent, as well
we use the edge-grafting transformations to identify the $n$-vertex trees of given bipartition having the second and third smallest Laplacian permanent. Similarly, we also characterize the $n$-vertex bipartite unicyclic graphs of given bipartition having the first, second and third smallest Laplacian permanent. Consequently, we identify the $n$-vertex bipartite unicyclic graphs with the first,
second and third smallest Laplacian permanent.

\section{\normalsize Three edge-grafting theorems on Laplacian permanent}\setcounter{equation}{0}
In this section, we introduce three edge-grafting transformations. We also study the property for each of the
three edge-grafting transformations. 
\begin{defi}
Let $uv$ be a pendant edge of an $n$-vertex bipartite graph $G$ with $d(u)=1, \, n\ge 3$. Let $w \, (\not=v)$ be a vertex of $G$
with $d(w)\geq d(v)$. Let $G [v\rightarrow w;1]$ be the graph
obtained from $G$ by deleting the edge $uv$ and adding the edge
$uw$. In notation,
$$
G[v \rightarrow w;1]=G-uv+uw \ \ 
$$
and we say $G[v \rightarrow w;1]$ is obtained from $G$ by
{\bf Operation I}.
\end{defi}
\begin{thm}\label{lem2.1}
Let $G$ and $G[v \rightarrow w;1]$ be the bipartite graphs defined as above. Then
$
  \per  L(G) > \per  L(G[v \rightarrow w;1]).
$
\end{thm}
\begin{proof}
Let $d_{G}(v)=r$ and $d_{G}(w)=t$. First we consider that $vw \not\in E_G$.
With an appropriate ordering of the vertices of $G$ and $G[v \rightarrow w;1]$ as $u,v,w,\ldots,$  we see that
\begin{equation*}
\begin{array}{cc}
  L(G)= &
\left(
 \begin{array}{cccc}
  1 & -1 & 0& {\bf 0} \\
  -1 & r & 0 & {\bf x}_1 \\
  0 & 0 & t & {\bf x}_2 \\
  {\bf 0} & {\bf y}_1 & {\bf y}_2 & \bf{A}
 \end{array}
\right)
\end{array}
\end{equation*}
and
\begin{equation*}
\begin{array}{cc}
  L(G[v\rightarrow w;1])= &
\left(
 \begin{array}{cccc}
  1 & 0 & -1& {\bf 0} \\
  0 & r-1 & 0 & {\bf x}_1 \\
  -1& 0 & t+1 & {\bf x}_2 \\
  {\bf 0} & {\bf y}_1 &{\bf y}_2 & \bf{A}
 \end{array}
\right).
\end{array}
\end{equation*}

Let $M_1$ (resp. $M_2$) be the matrix obtained from $L(G)$ (resp.
$L(G[v \rightarrow w;1])$) by eliminating the first row and the first
column. Let $N_1$ (resp. $N_2$) be the matrix obtained from $L(G)$
(resp. $L(G[v \rightarrow w;1])$) by eliminating the
first 2 rows and the first 2 columns. And let $N_2'$
be the matrix obtained from $M_2$ by eliminating its second row and
second column. Then we have
$$ \per L(G)= \per M_1+\per N_1,\ \ \ \ \per
L(G[v\rightarrow w;1])= \per M_2+\per N_2'.
$$

Set
$$
S_1=\per
\left(\begin{array}{cc}
        0 & {\bf x}_2 \\
        {\bf y}_2 & {\bf A}
      \end{array}
 \right), \ \ \ \ \ \ S_2=\per
\left(\begin{array}{cc}
        0 & {\bf x}_1 \\
        {\bf y}_1 & {\bf A}
      \end{array}
 \right).
$$
Note that
$$\per N_1=t\cdot \per A + S_1, \ \ \ \ \per N_2'=(r-1)\cdot \per A + S_2, \ \ \ \ \per M_2=\per M_1+\per N_2'-\per N_1.$$
Hence,
\begin{equation}\label{eq:2.1}
\per L(G)-\per L(G[v \rightarrow w;1])
=2((t-r)\cdot \per A +\per A +S_1 -S_2).
\end{equation}
By the choice of $S_1$ and $S_2$, we have
\begin{equation}\label{eq:2.2}
\per A+S_1 >S_2.
\end{equation}
Note that $t\geq r$,  by (\ref{eq:2.1}) and (\ref{eq:2.2}) we get
$
  \per L(G)-\per L(G[v \rightarrow w;1])>0.
$

Now consider $vw\in E_G.$ By an similar argument as in the proof of
the case $vw\not\in E_G,$ we can also get $\per L(G)-\per L(G[v \rightarrow
w;1])>0.$ We omit the procedure here.

This completes  the proof.
\end{proof}
\begin{defi}
Let $vw$ be an edge of a bipartite graph $U$ with
$d(w)\geq 2$. $G$ is obtained from $U$ and the star $S_{k+2}$ by
identifying $v$ with a pendant vertex of $S_{k+2}$ whose center is $u$.
Let $G[u\rightarrow w;2]$ be the graph obtained from $G$ by deleting all edges $uz, z\in W$ and adding all
edges $wz, z\in W$, where $W=N_{G}(u)\backslash\{v\}.$ In notation,
$$
G[u\rightarrow w;2]=G-\{uz:z\in W\}+\{wz:z\in W\}
$$
and we say $G[u\rightarrow w;2]$ is obtained from $G$ by
{\bf Operation I\!I}. Graphs $G$ and $G[u\rightarrow w;2]$ are depicted in Fig. 2.
\end{defi}
\begin{figure}[h!]
\begin{center}
\psfrag{a}{$w$}\psfrag{b}{$v$} \psfrag{c}{$u$}\psfrag{d}{$u_1$}\psfrag{u}{$U$}
\psfrag{e}{$u_k$} \psfrag{g}{$G$}\psfrag{h}{$G[u\rightarrow w;2]$}
\psfrag{f}{$u_2$} \psfrag{i}{$u_2$}
\includegraphics[width=110mm]{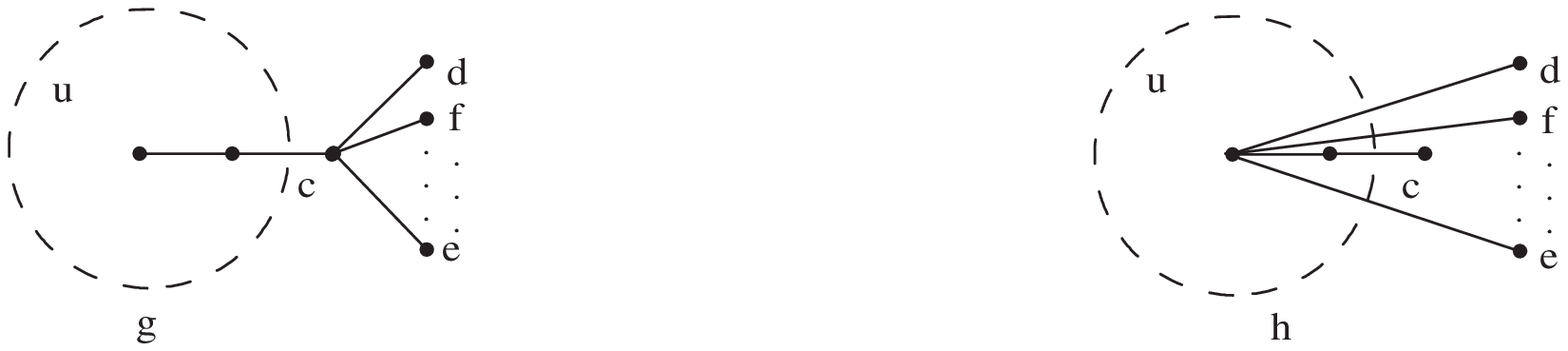}\\
\caption{$G\Rightarrow G[u\rightarrow w;2]$ by Operation I\!I.}
\end{center}
\end{figure}
\begin{thm}\label{lem2.2}
Let $G$ and $G[u\rightarrow w;2]$ be the bipartite graphs described
as above. Then $\per L(G) > \per L(G[u\rightarrow w;2]).$
\end{thm}
\begin{proof}
Let $N_G(u)\backslash\{v\}=\{u_1,u_2,\ldots,u_k\}$, where $k\ge 1,\,
u_1, u_2, \ldots, u_k$ are pendant vertices. With an appropriate
ordering of the vertices of $G$ and $G[u\rightarrow w;2]$ as
$u_1,u_2,\ldots,u_k,u,v,w,\ldots,$ we have
\begin{equation*}
\begin{array}{cc}
  L(G)= &
\left(
  \begin{array}{ccc|c|ccc}
    1 &   &   & -1 &   &   &   \\
      & \ddots &  &\vdots   &  &   &   \\
      &   & 1 & -1 &   &   &   \\
      \hline
    -1 & \cdots & -1 & k+1 & -1 &   &   \\
    \hline
      &  &   & -1 & d_v & -1 & {\bf x}_1 \\
      &   &   &  & -1 & d_w & {\bf x}_2 \\
      &   &   &  & {\bf y}_1 & {\bf y}_2 & \textbf{A }
  \end{array}
\right)
\end{array}
\end{equation*}
and
\begin{equation*}
\begin{array}{cc}
  L(G[u\rightarrow w;2])= &
\left(
  \begin{array}{ccc|cc|cc}
    1 &   &   &   &   & -1 &   \\
      & \ddots &   &   &   & \vdots &   \\
      &   & 1 &   &   & -1 &   \\
    \hline
      &   &   & 1 & -1 & 0 &   \\
      &   &   & -1 & d_v & -1 & {\bf x}_1 \\
    \hline
    -1 & \cdots & -1 & 0 & -1 & d_w+k & {\bf x}_2 \\
      &   &   &   &{\bf y}_1 & {\bf y}_2 & \bf{A}
  \end{array}
\right),
\end{array}
\end{equation*}
where $d_v$ and $d_w$ are degrees of $v$ and $w$ in $G$ with $d_w\ge 2$.

Let $D_1$ (resp. $D_2$) be the matrix obtained from $L(G)$
(resp. $L(G[u\rightarrow w;2])$) by eliminating the first $k+1$ rows and
the first $k+1$ columns; let $M_1$ (resp. $M_2$) be the matrix
obtained from $D_1$ (resp. $D_2$) by eliminating the
first row and the first column. Let $M$ be the matrix
obtained from $D_1$ by eliminating the second row and the second column, and let $M(i)$ (resp. $N(i))\, (1\le i\le k)$  be the
matrix obtained from $L(G[u\rightarrow w;2])$ by eliminating the rows and columns corresponding to $u_1, u_2, \ldots, u_{k-i}$ and $u$ (resp. $u_1, u_2, \ldots, u_{k-i}, u$ and $v$).

It is routine to check that $\{\per M(i), 1\leq i \leq k \}$ and $\{\per N(i), 1\leq i \leq k \}$, respectively, have the recurrence relation
$$
\per M(i)= \per M(i-1)+\per M,  \ \ \ \ \per N(i)= \per N(i-1)+\per A, \ \ \ \ 2\leq i \leq k
$$
with initial value $\per M(1)=\per D_2+\per M$ and $\per N(1)=\per M_2+\per A$. Hence, we have
\begin{equation}\label{eq:2.4}
\per M(k)=\per D_2+k\cdot \per M,\ \ \ \ \ \per N(k)=\per M_2+k \cdot\per A.
\end{equation}
By expanding the permanent of $L(G)$ along the first $(k+1)$ rows we
obtain
\begin{equation}\label{eq:2.6}
\per L(G)=(2k+1)\cdot\per D_1+\per M_1,
\end{equation}
and by expanding the permanent of $L(G[u\rightarrow w;2])$ along the row
corresponding to $u$, we get
$$
\per L(G[u\rightarrow w;2])=\per M(k)+\per N(k).
$$
Note that
$\per D_2=\per D_1+k\cdot \per M$ and $\per M_2=\per M_1+k\cdot \per A$. Hence, by (\ref{eq:2.4}) we have
\begin{equation}\label{eq:2.7}
\per L(G[u\rightarrow w;2])=\per D_1+2k\cdot \per M+\per M_1+2k\cdot \per A.
\end{equation}
For convenience, denote by $D_1'$ the matrix obtained from $D_1$
by replacing $d_w$ with $d_w-1$. In view of (\ref{eq:2.6}) and
(\ref{eq:2.7}), we get
\begin{eqnarray*}
  \per L(G)-\per L(G[u\rightarrow w;2])&=&2k\cdot\per D_1-2k\cdot\per M-2k\cdot\per A   \\
   &=&2k(\per D_1'+\per M)-2k\cdot\per M-2k\cdot\per A  \\
   &=&2k(\per D_1'-\per A)  \\
   &\geq&2k[(d_v(d_w-1)+1)\cdot\per A-\per A]\\
   &>&0.
\end{eqnarray*}

This completes the proof.
\end{proof}
\begin{defi}
Let $G$ be an $n$-vertex graph obtained from $C_{2k}=v_1v_2\ldots
v_i\ldots v_j\ldots v_{2k}v_1\, (k \ge 3)$ and two stars $S_{n_i+1},
S_{n_j+1}$ by identifying $v_i$ (resp. $v_j$) with the center of
$S_{n_i+1}$ (resp. $S_{n_j+1}$), where $4<i<j,\, n=2k+n_i+n_j;$\,
see Fig. 3. Let $G'=G- v_1v_2+v_1v_4$. Then we say that $G'$ is
obtained from $G$ by {\bf Operation I\!I\!I}.
\end{defi}
\begin{figure}[h!]
\begin{center}
\psfrag{a}{$v_1$}\psfrag{b}{$v_2$}
\psfrag{c}{$v_3$}\psfrag{d}{$v_4$}
\psfrag{p}{$C_{2k}$}\psfrag{q}{$C_{2k-2}$}
\psfrag{e}{$G$}\psfrag{f}{$G'$}
\psfrag{g}{$n_j$}\psfrag{h}{$n_i$}
\psfrag{l}{$v_i$}\psfrag{m}{$v_j$}
\includegraphics[width=100mm]{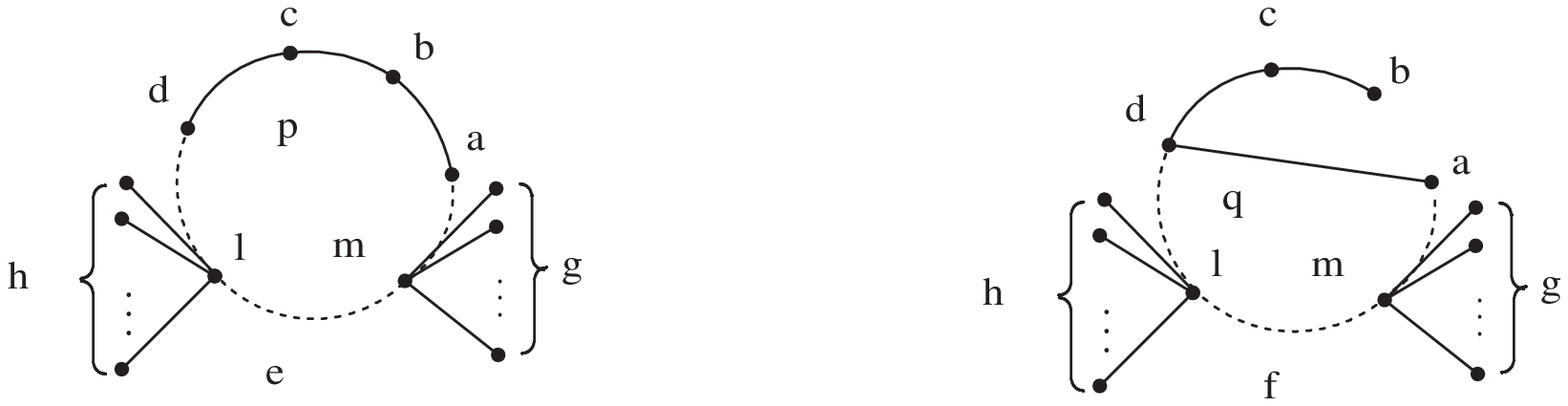}\\
\caption{$G\Rightarrow G'$ by Operation I\!I\!I.}
\end{center}
\end{figure}
\begin{thm}\label{lem2.3}
Let $G$ and $G'$ be the bipartite unicyclic graphs described as above.
Then $\per L(G) > \per L(G').$
\end{thm}
\begin{proof}
Ordering the vertices of $G$ as
$v_1,v_2,v_3,v_4,\ldots$, we have
\begin{equation*}
\begin{array}{cc}
  L(G)= &
\left(
  \begin{array}{ccccc}
    2 & -1 &   &   & {\bf x}_1 \\
    -1 & 2 & -1 &   &   \\
      & -1 & 2 & -1 &   \\
      &   & -1 & 2 & {\bf x}_4 \\
    {\bf y}_1 &   &   & {\bf y}_4 & {\bf A}
  \end{array}
\right).
\end{array}
\end{equation*}

Ordering the vertices of $G'$ as
$v_2,v_3,v_4,v_1,\ldots$, we see that
\begin{equation*}
\begin{array}{cc}
  L(G')= &
\left(
  \begin{array}{ccccc}
    1 & -1 &   &   &   \\
    -1 & 2 & -1 &   &   \\
     & -1 & 3 & -1 & {\bf x}_4 \\
      &   & -1 & 2 & {\bf x}_1 \\
      &   &  {\bf y}_4 & {\bf y}_1 & {\bf A}
  \end{array}
\right).
\end{array}
\end{equation*}

Let $D$ (resp. $D'$) be the matrix obtained from $L(G)$ (resp.
$L(G')$) by eliminating the rows and columns corresponding to $v_2$
and $v_3$; and let $M_1$ (resp. $M_2$) be obtained from $L(G)$
(resp. $L(G')$) by eliminating rows and columns corresponding to
$v_2, v_3, v_4$ (resp. $v_1, v_2, v_3$). And for convenience, denote
\begin{equation*}
\begin{array}{ccc}
 \begin{array}{cc}
  N_1= &
\left(
  \begin{array}{cc}
    0 & {\bf x}_1 \\
    {\bf y}_4 &\bf{ A} \\
  \end{array}
\right),
\end{array}
 & &
\begin{array}{cc}
  N_2= &
\left(
  \begin{array}{cc}
    0 & {\bf x}_4 \\
   {\bf y}_1 &\bf{ A} \\
  \end{array}
\right).
\end{array}
\end{array}
\end{equation*}

Expanding the rows corresponding to $v_2$ and $v_3$ of $L(G)$
and $L(G')$, respectively, yields
\begin{eqnarray*}
\per L(G)&=&-\per N_1+2\per M_2+\per A+5\per D+2\per M_1-\per N_2,\\
\per L(G')&=&3\per D'+\per M_1.
\end{eqnarray*}
Together with $\per D'=\per D+\per M_1+\per A-\per N_1+\per N_2
$, we obtain
\begin{equation*}
\per L(G')=3 \cdot\per D +4 \cdot\per M_1+ 3\cdot\per A -3\cdot\per
N_1 -3\cdot\per N_2
\end{equation*}
and hence
\begin{align}\label{eq:2.8}
  \per L(G)-\per L(G') =&
  2\per D+2\per M_2+2\per N_1+2\per N_2-2\per M_1-2\per A \notag\\
   =&2(\per D-\per M_1)+2(\per M_2-\per A)+2\per N_1+2\per N_2.
\end{align}
By ordering the vertices of $G$ as $v_1,v_2,v_3,v_4,\ldots, v_i,
\ldots, v_j, \ldots, v_{2k},\ldots,$ and by direct calculation, we
have $\per  N_1=\per N_2=-1.$ And note that
$$
  (\per D-\per M_1)+(\per M_2-\per A) = \per
  \left(
    \begin{array}{ccc}
      3 & 0 & {\bf x}_1 \\
      0 & 1 & {\bf x}_4 \\
      {\bf y}_1 & {\bf y}_4 & \bf{A} \\
    \end{array}
  \right)\\
   \geq 3.
$$
Together with (\ref{eq:2.8}), we have
$$
 \per L(G)-\per L(G')\geq 2\times 3+2(-1)+2(-1)=2>0.
$$
This completes the proof.
\end{proof}

\section{\normalsize Applications}
\subsection{\normalsize Laplacian permanents of trees among $\mathscr{T}_n^{p,q}$}

\setcounter{equation}{0}
We denote by $D(p,q)$ a \textit{double star} with $n$ vertices,
which is obtained from an edge $vw$ by attaching $p-1$ (resp. $q-1$)
pendant vertices to $v$ (resp. $w$), where $n=p+q$. Let $D'(p-1,q-1)$
(resp. $D''(p-1,q-1)$) be an $n$-vertex tree obtained from
$D(p-1,q-1)$ by attaching a pendant path of length 2 to $w$ (resp.
$v$). Graphs $D(p,q),\, D'(p-1,q-1),\, D''(p-1,q-1)$ are depicted in Fig. 4. Let $T(n,k,a)$ be an $n$-vertex tree obtained by
attaching $a$ and $n-k-a$ pendant vertices to the two end-vertices of $P_k$, respectively. In particular, $D(p,q)=T(n,2,p-1).$
\begin{figure}[h!]
\begin{center}
\psfrag{a}{$p-1$}\psfrag{b}{$q-1$}\psfrag{c}{$p-2$}\psfrag{d}{$q-2$}\psfrag{e}{$q-2$}
\psfrag{v}{$v$}\psfrag{w}{$w$}
\psfrag{A}{$D(p, q)$}\psfrag{C}{$D'(p-1,q-1)$}\psfrag{B}{$D''(p-1,q-1)$}
  \includegraphics[width=119mm]{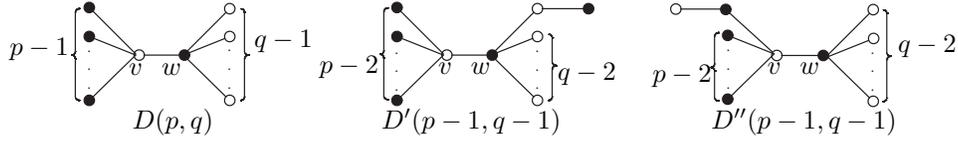}\\
  \caption{Trees $D(p, q), D'(p-1,q-1)$ and $D''(p-1,q-1).$ }
\end{center}
\end{figure}

The following lemma is routine to check.
\begin{lem}
Let $p$ and $q$ be positive integers, then
$$
\per L(D'(p-1,q-1))=(2p-3)(6q-5)+3, \,
\per L(D''(p-1,q-1))=(2q-3)(6p-5)+3.
$$
\end{lem}

From Lemma 3.1, a direct calculation yields
\begin{equation}\label{a}
\text{$\per L(D'(p-1,q-1))< \per L(D''(p-1,q-1))$\ \ \ \ \  for $q>p>2.$}
\end{equation}

From \cite{B-Z1} we know that $D(p,q)$ minimizes the Laplacian
permanent of trees among $\mathscr{T}_n^{p,q}$. In this subsection, we use a new and unified method to determine
the tree in $\mathscr{T}_n^{p,q}$ which has the first, second, and third smallest Laplacian permanent, respectively.
\begin{thm}[\cite{B-Z1}]
Let $T$ be a tree with a $(p,\, q)$-bipartition. Then
$$
\per L(T) \geq (2p-1)(2q-1)+1
$$
with equality if and only if $T$ is a double-star $D(p,q)$.
\end{thm}
\begin{proof}
Choose a tree $T$ with a $(p,q)$-bipartition such that its Laplacian permanent
is as small as possible. Let $V_1$, $V_2$ be the bipartition of the
vertices of $T$ with $V_1=\{v_0,v_1,\ldots,v_{p-1}\}$, $V_2=\{ u_0,u_1,\ldots,u_{q-1}\}$. For
convenience, let $v_0$ (resp. $u_0$) be the vertex of maximal degree
among $V_1$ (resp. $V_2$) in $T$ and let $A=N_T(v_0)\cap PV(T)$. 

Hence, in order to complete the
proof, it suffices to show the following claims.
\begin{figure}[h!]
\begin{center}
\psfrag{a}{$v_0$} \psfrag{b}{$u_1$} \psfrag{c}{$v_1$}
\psfrag{d}{$u_{r-1}$} \psfrag{e}{$v_{r-1}$} \psfrag{f}{$u_0$}
\psfrag{g}{$T(n,2r,s)$}\psfrag{m}{$s$} \psfrag{n}{$n-2r-s$}
  \includegraphics[width=70mm]{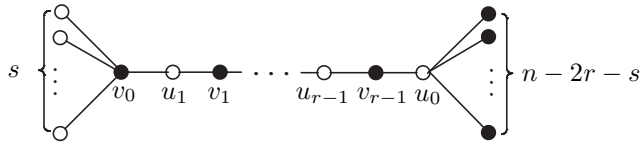}\\
  \caption{Tree $T(n,2r,s)$ with some labelled vertices. }
\end{center}
\end{figure}
\begin{claim}
$T \cong T(n,2r,s)$ (see Fig. 5) with $r\geq 1$ and $s \geq 0$.
\end{claim}

\noindent{\bf Proof of Claim 1.} Otherwise, $T$ must contain a
pendant vertex $w \not\in N_T(u_0)\cup N_T(v_0)$. Without loss of
generality, we may assume $w \in V_2$ and its unique neighbor is
$w'$. Using Operation I, let $ T_0 = T- ww'+ wv_0. $ By Theorem 2.2
we have $\per L(T_0)< \per L(T)$, a contradiction to the choice of
$T$.

This completes the proof of Claim 1. \qed
\begin{claim}
In the tree described as above, $u_0$ is adjacent to $v_0$.
\end{claim}

\noindent{\bf Proof of Claim 2.} If not, then $d(u_0,v_0)\geq
3$. Note that $v_0$ is of the maximal degree vertex, hence $d_T(v_0)\not =1$, which implies $A\not=\emptyset.$ Using Operation I\!I, let
$$
T_1= T-\{ v_0z:z\in A\}+ \{ v_1z:z\in A\}
$$
then $\per L(T_1)< \per L(T)$ by Theorem 2.4, which contradicts the
choice of $T$. This completes the proof of Claim 2. \qed

By Claims 1 and 2, we get that $T\cong D(p,q)$. By direct computing,
we have
$$
\per L(D(p,q))=(2p-1)(2q-1)+1.
$$

This completes the proof of Theorem 3.2.
\end{proof}

\setcounter{claim}{0}

\begin{thm}
Among $\mathscr{T}_n^{p,q}$.
\begin{wst}
  \item[{\rm (i)}] If $p =2$, then all the members in $\mathscr{T}_n^{2,n-2}$ are ordered as follows:
\begin{equation}\label{eq:3.2}
\begin{split}
&\per L(D(2,n-2))=\per L(T(n,3,0))< \per L(T(n,3,1))< \per L(T(n,3,2)\notag\\
&\ \ \   < \cdots <\per L(T(n,3,i))<\cdots <\per
L(T(n,3,\lfloor \frac{ n-3 }{2}\rfloor)).\\
\end{split}
\end{equation}
  \item[{\rm (ii)}] If $p >2$,
\begin{wst}
  \item[{\rm (a)}] for $T \in \mathscr{T}_n^{p,q}\setminus\{D(p,q)\}$, we have
$ \per L(T)\geq  (2p-3)(6q-5)+3.$ The equality holds if and only if
$T \cong D'(p-1,q-1).$
  \item[{\rm (b)}] for $T \in \mathscr{T}_n^{p,q}\setminus\{D(p,q), D'(p-1,q-1)\}$ with $q>p$, we have
$ \per L(T)\geq (2q-3)(6p-5)+3. $ The equality holds if and only if
$T \cong D''(p-1,q-1).$
\end{wst}
\end{wst}
\end{thm}
\begin{proof} (i)\ If $p=2$, then
$$
\mathscr{T}_n^{2,n-2}=\left\{T(n,3,0),T(n,3,1),\ldots,
T(n,3,i),\ldots,T(n,3,\left\lfloor\frac{n-3
}{2}\right\rfloor)\right\}.
$$
By a simple calculation, we get
\begin{eqnarray*}
  \per L(T(n,3,i)) =-8(i- \frac{ n-3 }{2})^2+2(n-3)^2+6n-14.
\end{eqnarray*}

Consider the function $f(x)=-8(x- \frac{ n-3 }{2})^2+2(n-3)^2+6n-14$
in $x$ with $0 \leq x \leq \lfloor \frac{n-3}{2} \rfloor.$ By the
monotonicity of $f(x)$, we have
\[\label{eq:3.3}
f(0)< f(1)< \cdots < f(i-1)< f(i)<f(i+1)<\cdots<f(\lfloor\frac{
n-3}{2}\rfloor).
\]
Note that $T(n,3,0)\cong D(2,n-2)$, together with (\ref{eq:3.3}) we get that
(i) holds.

(ii)\ Choose $T\in \mathscr{T}_n^{p,q}\setminus \{D(p,q)\}$ such
that its Laplacian permanent is as small as possible. Let $V_1$,
$V_2$ be the bipartition of the vertices of $T$ with $V_1=\{ v_0,v_1,\ldots,v_{p-1}\}$, $V_2=\{
u_0,u_1,\ldots,u_{q-1}\}$. For convenience, let $v_0$ (resp. $u_0$)
be the vertex of maximal degree among $V_1$ (resp. $V_2$) in $T$ and
let $A=N_T(v_0)\cap PV(T)$. In order to
complete the proof, it suffices to show the following claims.
\begin{claim}
$u_0v_0\in E_T.$
\end{claim}

\noindent{\bf Proof of Claim 1.}\ If not, then $d_T(u_0,v_0)\ge 3$.
In this case, we are to show that $T\cong T(n,2r,s)$ (see Fig. 5).
Otherwise, $T$ must contain a pendant vertex $w \not\in N_T(u_0)\cup
N_T(v_0)$. Assume that the unique neighbor of $w$ is $w'$. Using
Operation I, let $ T' = T- ww'+ wv_0 $ if $w\in V_2$ and $ T' = T-
ww'+ wu_0 $ otherwise. Note that $T'\in \mathscr{T}_n^{p,q}\setminus
\{D(p,q)\}$, by Theorem 2.2 we have $\per L(T')< \per L(T)$, a
contradiction to the choice of $T$. Hence, $T\cong T(n,2r,s)$. On
the one hand, $d_T(u_0,v_0)\ge 3$, hence $r\ge 2$; on the other
hand, $v_0$ is of the maximal degree vertex in $V_1$ of $T$, hence
$s\ge 1.$ Therefore, $A\not=\emptyset.$ Using Operation I\!I, let
$$
T_0 = T- \{v_0z: z \in A \} + \{v_1z: z \in A\}.
$$
We also have $T_0\in \mathscr{T}_n^{p,q}\setminus \{D(p,q)\}$. In
view of Theorem 2.4, $\per L(T_0)< \per L(T)$, which also
contradicts the choice of $T$.

This completes the proof of Claim 1. \qed
\begin{claim}In the tree $T$ described as above, there exists a pendant vertex, say $w$, in $V_T$ such that
$d_T(w, u_0)=2, d_T(w, v_0)=3$ or $d_T(w, v_0)=2, d_T(w, u_0)=3.$
Furthermore, for all $v\in PV(T)\setminus\{w\},\, v$ is adjacent to
either $u_0$ or $v_0$.
\end{claim}

\noindent{\bf Proof of Claim 2.}\ Note that $T\not\cong D(p,q)$,
hence there must exist a vertex (not necessary a pendant vertex),
say $w$, such that $d_T(w, u_0)=2, d_T(w, v_0)=3$ or $d_T(w, v_0)=2,
d_T(w, u_0)=3.$ With loss of generality, we assume that $d_T(w,
u_0)=2, d_T(w, v_0)=3$. If $w$ is not a pendant vertex, then $w$ is
on a path which joins $u_0$ and a pendant vertex, say $r$. Denote
the unique neighbor of $r$ by $r'$. Let $T'=T-rr'+ru_0$ if $r\in
V_1$ and $T'=T-rr'+rv_0$ otherwise. It is routine to check that
$T'\in \mathscr{T}_n^{p,q}\setminus \{D(p,q)\}$. By Theorem 2.2,
$\per L(T')<\per L(T),$ a contradiction to the choice of $T$. Hence,
$w$ must be a pendant vertex.

In what follows, we should show that for all $v\in
PV(T)\setminus\{w\}$, either $vu_0\in E_T$ or $vv_0\in E_T.$ In
fact, if there exist a vertex, say $y$, in $PV(T)\setminus\{w\}$
such that $yu_0, yv_0\not\in E_T$. Denote the unique neighbor of $y$
by $y'$. Let $ \hat{T}=T-yy'+yu_0$ if $y\in V_1$ and
$\hat{T}=T-yy'+yv_0$ otherwise. It is straightforward to check that
$\hat{T}\in \mathscr{T}_n^{p,q}\setminus \{D(p,q)\}$. By Theorem
2.2, $\per L(\hat{T})<\per L(T),$ a contradiction to the choice of
$T$.

This completes the proof of Claim 2. \qed

By Claims 1 and 2 we obtain
$
T\cong D'(p-1,q-1),$  or $T\cong D''(p-1,q-1).$ If $p=q$, then
$D'(p-1,q-1)\cong D''(p-1,q-1).$ Together with Lemma 3.1, (ii) holds obviously in this case. If $p<q$, then
combining with Lemma 3.1 and Inequality (3.1), (ii) follows
immediately.

This completes the proof.
\end{proof}

\noindent{\bf Remark 1.} We know from \cite{W-G-Y} that $D(p,q)
\prec D'(p-1,q-1)\prec T$ for all $T\in
\mathscr{T}_n^{p,q}\setminus\{D(p,q), D'(p-1,q-1)\}$. In view of
Theorems 3.2 and 3.3, $D(p,q)$ (resp. $D'(p-1,q-1)$) is the tree
with a $(p, q)$ bipartition which has the smallest (resp. second
smallest) Laplacian permanent. Hence, our result support the
conjecture that trees minimizing the Laplacian permanent usually
simultaneously minimize the Laplacian coefficients, and vice versa.
Furthermore, in view of Theorems 3.2 and 3.3 it is natural to
conjecture that $D(p,q) \prec D'(p-1,q-1)\prec D''(p-1,q-1)\prec T$
for all $T\in \mathscr{T}_n^{p,q}\setminus\{D(p,q), D'(p-1,q-1),
D''(p-1,q-1)\}$ with $q>p$.

\subsection{\normalsize Laplacian permanent of trees with diameter at least $d$}
\setcounter{claim}{0}
Let $Q_{n}$ be the matrix obtained from
$L(P_{n+1})$ by eliminating row 1 and column 1. It is routine to check that $\per Q_1=1$
and $\per Q_2=3$. In particular, define $\per Q_0=1$. We know \cite{B-Z1} that
\[\label{eq:3.5}
\per Q_{n}=\frac{1}{2}\left(1+\sqrt{2}\right)^{n}+\frac{1}{2}\left(1-\sqrt{2}\right)^{n}
\]
and
\[\label{eq:3.6}
\per
L(P_{n})=\frac{2-\sqrt{2}}{2}\left(1+\sqrt{2}\right)^{n}+\frac{2+\sqrt{2}}{2}\left(1-\sqrt{2}\right)^{n}.
\]
\begin{lem}[\cite{B-Z1}]\label{lem3.4}
Let $n,j$ and $k$ be positive integers with $1\leq k< j\leq
\frac{1}{2}(n+1)$. Then
$
(-1)^{k}(\per Q_{j-1}\per Q_{n-j}- \per Q_{k-1}\per Q_{n-k})>0.
$
\end{lem}
\begin{lem}[\cite{G-H-L}]\label{lem3.5}
Let $uv$ be the only non-pendant edge incidence with $v$ in a tree
$T$ and let $A=N_T(v) \setminus \{ u\}$. Let
$
T'=T-\{ vz:z\in A\}+\{ uz:z\in A\},
$
then we have
$
\per L(T')< \per L(T).
$
\end{lem}

In this subsection, we use a new method to prove the following known result.
\begin{thm}[\cite{B-Z1}]
Let $d$ be a positive integer, and let $T$ be a tree with $n$
vertices having diameter at least $d$. Then
$$
 \per L(T)\geq \left(n-d+\frac{\sqrt{2}}{2}\right)\left(1+\sqrt{2}\right)^{d-1}+\left(n-d-\frac{\sqrt{2}}{2}\right)\left(1-\sqrt{2}\right)^{d-1}.
$$
The equality holds if and only if $T\cong T_{n,d,2};$ see Fig. 1.
\end{thm}
\begin{proof}
Choose an $n$-vertex tree $T$ of diameter at least $d$ such that its
Laplacian permanent is as small as possible. If $T\cong P_{d+1}$,
then our result holds by Theorem 1.1. Hence, in what follows we
consider that $T\not\cong P_{d+1}$. If $T$ contains just two pendant
vertices, i.e., $T\cong P_{n}=v_1v_2\ldots v_i\ldots v_n$. Let
$T'=T-v_1v_2+v_iv_1$. Obviously, $T'$ is of diameter at least $d$.
By Theorem 2.2, we have $\per L(T')<\per L(T)$, a contradiction to
the choice of $T.$ Hence, $T$ contains at least 3 pendant vertices.
That is to say, the maximal vertex degree in $T$ is of at least 3.
Without loss of generality, we may assume that $w$ is just of the
maximal degree vertex. Let $P'=v_1v_2\ldots v_i\ldots v_{l+1}$ be
one of the longest path contained in $T$, where $l\ge d$. In order
to complete the proof, it suffices to show the following claims.

\begin{claim}
$T\cong T_{n,l,i}$, where $T_{n,l,i}$ is obtained from $P'$ by
inserting $n-l-1$ pendant vertices at $v_i$, $i\in \{2,3,\ldots,
\lfloor(l+2)/2\rfloor\}.$
\end{claim}

\noindent{\bf Proof of Claim 1.}\ First we show that all the pendant
vertices excluding the endvertices of $P'$ are adjacent to $w$.
Assume to the contrary that $v\in PV(T)\setminus \{v_1, v_{l+1}\}$
satisfying $vw\not \in E_T.$ Denote the unique neighbor of $v$ by
$v'$. Set $T'=T-vv'+wv$. It is straightforward to check that $T'$ is
of an $n$-vertex tree of diameter at least $d$. By Theorem 2.2, we
have $\per L(T')<\per L(T)$, a contradiction to the choice of $T.$

Now we show that $w$ is on the path $P'.$ Assume that $w$ is not on
the path $P'$, then $T$ must be the tree obtained by joining the
center of a star $S$ and a vertex of $P'$ by a path of length at
least 1. Denote the unique neighbor of $w$ which is not a pendant
vertex by $w'$. Set $A=N_T(w) \setminus \{w'\}$. Let $ T'=T-\{
wz:z\in A\}+\{ w'z:z\in A\}. $ It is easy to see that $T'$ is a tree
of diameter at least $d$. By Lemma 3.5, $\per L(T')<\per L(T)$, a
contradiction.

This completes the proof of Claim 1. \qed

\begin{claim}
In the tree $T_{n,l,i}$ described as above, we have $l=d, i=2$, i.e., $T_{n,l,i}\cong T_{n,d,2}$.
\end{claim}

\noindent{\bf Proof of Claim 2.}\ If not, then $l\geq d+1$. In the
tree described above, let $ T_2=T_{n,l,i}-v_{l+1}v_l+v_iv_{l+1}. $
It is easy to see that $T_2$ is an $n$-vertex tree of diameter at
least $d$. By Theorem 2.2, we have $\per L(T_2)<\per L(T_{n,l,i})$,
a contradiction to the choice of $T_{n,l,i}$. So we obtain $T\cong
T_{n,d,i}$.

Expanding the permanent of
$L(T_{n,d,i})$ along the row corresponding to vertex $v_i$ gives
$$
\per L(T_{n,d,i})= \per L(P_{d+1}) +2 (n-d-1)\per Q_{i-1}\per Q_{d-i+1}.
$$
This gives
\begin{eqnarray*}
\per L(T_{n,d,j})-\per L(T_{n,d,2})& = &2(n-d-1)(\per Q_{j-1}\per
Q_{d-j+1}
- \per Q_{2-1}\per Q_{d-2+1})\\
&>&0
\end{eqnarray*}
for $j=3,4,\ldots, \lfloor \frac{1}{2}(d+2)\rfloor$ and the last
inequality follows by Lemma 3.4. \qed

In view of (\ref{eq:3.5}) and (\ref{eq:3.6}), we have
\[\label{eq:3.7}
\per
L(T_{n,d,2})=(n-d+\frac{\sqrt{2}}{2})(1+\sqrt{2})^{d-1}+(n-d-\frac{\sqrt{2}}{2})(1-\sqrt{2})^{d-1}.
\]

By Claims 1 and 2 and Eq. (\ref{eq:3.7}), Theorem 3.6 follows immediately.
\end{proof}

\subsection{\normalsize Lower bounds for the Laplacian permanent of graphs in $\mathscr{U}_n^{p,q}$}

\setcounter{claim}{0} In this subsection, we are to determine sharp
lower bounds for the Laplacian permanent of graphs in
$\mathscr{U}_n^{p,q}$. Let
$C_4(1^{s_1}k_1,\, 1^{s_2}k_2,\, 1^{s_3}k_3,\, 1^{s_4}k_4)$ be the
graph obtained from $C_4=v_1v_2v_3v_4v_1$ by inserting $s_i$ pendant vertices at $v_i$ and joining $v_i$ to the center of a star
$S_{k_i}$ by an edge, $i=1,2,3,4;$ Fig. 6. In
particular, let $B(p,\, q)=C_4(1^{p-2}0,\, 1^{q-2}0,\,1^00,\,1^00)$.
\begin{figure}[h!]\begin{center}
\psfrag{s}{$q-2$} \psfrag{t}{$p-2$}
\psfrag{a}{$v_1$}\psfrag{b}{$v_2$}
\psfrag{c}{$v_3$}\psfrag{d}{$v_4$}
\psfrag{e}{$s_1$}\psfrag{f}{$t_1$}
\psfrag{g}{$s_2$}\psfrag{h}{$t_2$}
\psfrag{i}{$s_3$}\psfrag{j}{$t_3$}
\psfrag{k}{$s_4$}\psfrag{l}{$t_4$}
\psfrag{m}{$k_1-1$}\psfrag{o}{$k_2-1$}
\psfrag{p}{$k_3-1$}\psfrag{q}{$k_4-1$}
\psfrag{A}{$C_4(1^{s_1}2^{t_1},1^{s_2}2^{t_2},1^{s_3}2^{t_3},1^{s_4}2^{t_4})$}
\psfrag{B}{$C_4'(1^{s_1}k_1,1^{s_2}k_2,1^{s_3}k_3,1^{s_4}k_4)$}
\psfrag{C}{$B(p,q)$}
  \includegraphics[width=60mm]{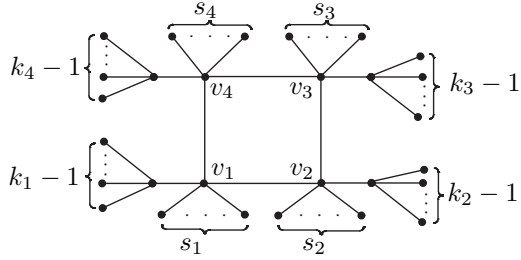}\\
  \caption{Graph $C_4(1^{s_1}k_1,1^{s_2}k_2,1^{s_3}k_3,1^{s_4}k_4).$}
\end{center}
\end{figure}
\begin{thm}
For any $G \in \mathscr{U}_n^{p,q},$ one has $ \per L(G)\geq
20(p-1)(q-1)+4n. $ The equality holds if and only if $G\cong
B(p,q).$
\end{thm}
\begin{proof}
Choose $G \in \mathscr{U}_n^{p,q}$ such that its Laplacian permanent
is as small as possible. If $n=4,\, 5$, $\mathscr{U}_4^{p,q}=\{
B(2,2)\}$ and $\mathscr{U}_5^{p,q}=\{ B(2,3)\}$, our result holds obviously. Hence in what
follows we consider $n\geq 6.$ Let $C_r$ be the unique cycle
contained in $G$. Note that, by Theorem 2.6, $G\not\cong C_n$.
Hence, $PV(G)\not=\emptyset.$ Let $V_1$, $V_2$ be the bipartition of
$V_G$ such that $|V_1|=p$ and $|V_2|=q$ with $v_0$ (resp. $u_0$)
being of the maximal degree vertex among $V_1$ (resp. $V_2$) in $G$.
\begin{claim}
For all $u\in PV(G)$, either $uu_0\in E_G$ or $uv_0\in E_G.$
\end{claim}

\noindent{\bf Proof of Claim 1.}\ If not, then there exists a
pendant vertex, say $u$, such that $u$ is not in $N_G(u_0)\cup
N_G(v_0)$. Denote the unique neighbor of $u$ by $u'$. Using
Operation I, let $G'=G-uu'+uu_0$ if $u\in V_1$ and $G'=G-uu'+uv_0$
otherwise. It is routine to check that $G' \in \mathscr{U}_n^{p,q}$.
By Theorem 2.2, $\per L(G')< \per L(G),$ a contradiction to the
choice of $G$.

This completes the proof of Claim 1. \qed

Let $d(u,C_r)=\min\{d(u,v): \ v\in V_{C_r}\}$. In particular, if $u$
is on $C_r$, then $d(u,C_r)=0$.
\begin{claim}
$d(v_0,C_r) = d(u_0,C_r) =0$.
\end{claim}

\noindent{\bf Proof of Claim 2.}\ Here we only show that $d(v_0,C_r)
=0$ by contradiction. With the same method, we can also show $d(u_0,C_r) =0$.

Assume that $d(v_0,C_r) =t \ge 1$. Set $A= PV(G)\cap N_G(v_0).$ By
Claim 1, we have $A\not=\emptyset.$ Let $P_t=v_0w_1w_2\ldots
w_{t-1}w_t$ be the shortest path connecting $v_0$ and the cycle
$C_r$, where $w_t$ is on $C_r$. Let $u\in V_{C_r}\cap N_G(w_t).$
Using Operation I\!I, let
$$
\bar{G}=\left\{
          \begin{array}{ll}
            G-\{ zv_0: z\in A\}+\{ zu: z\in A\}, & \hbox{if $t=1$;} \\
            G-\{ zv_0: z\in A\}+\{ zw_2: z\in A\}, & \hbox{if $t\ge 2$.}
          \end{array}
        \right.
$$
It is routine to check that $\bar{G}\in \mathscr{U}_n^{p,q}$. By
Theorem 2.4, we have $\per L(\bar{G})< \per L(G)$, a contradiction to the
choice of $G$. \qed

\begin{claim}
$r\leq 6$.
\end{claim}

\noindent{\bf Proof of Claim 3.}\ If not, then $r\geq 8$. By the
structure of $G$ described as above, then there must exist four
consecutive vertices, say $u_{k_1}, v_{k_1}, u_{k_2}, v_{k_2}$ on
the cycle $C_r$ such that $u_0,v_0\not\in \{u_{k_1}, v_{k_1},
u_{k_2}, v_{k_2}\}$. Without loss of generality assume that
$v_{k_1}, v_{k_2} \in V_1$ and $u_{k_1},u_{k_2}\in V_2$. Using
Operation I\!I\!I, let $ G_0=G- u_{k_1}v_{k_1}+u_{k_1}v_{k_2}. $ It
is routine to check that $G_0\in \mathscr{U}_n^{p,q}$. By Theorem
2.6 $\per L(G_0) < \per L(G),$ a contradiction to the choice of $G$.
This completes the proof of Claim 3.

Hence, by Claims 1-3, we obtain
\begin{itemize}
  \item $r=4$, then $G\cong B(p,q).$
  \item $r=6$, then $G\cong G_1$ or $G_2$, where $G_1, G_2$ are depicted in Fig. 7.
\end{itemize}
\begin{figure}[h!]
\begin{center}
\psfrag{a}{$v_0$} \psfrag{b}{$u_0$} \psfrag{c}{$q-3$}
\psfrag{d}{$p-3$} \psfrag{e}{$G_1$} \psfrag{f}{$G_2$}
\includegraphics[width=110mm]{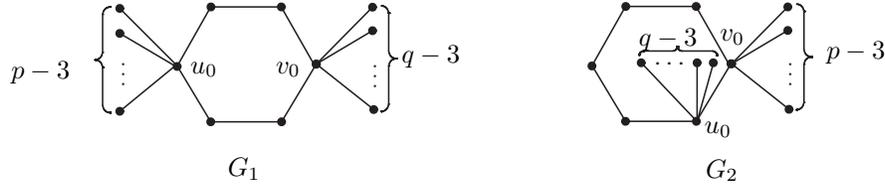}\\
\caption{Graphs $G_1$ and $G_2$.}
\end{center}
\end{figure}

If $G\cong G_2$, using Operation I\!I\!I on $G_2$, we obtain graph
$C_4(1^02,1^{q-3}0,1^{p-3}0,1^00)$ which is in
$\mathscr{U}_n^{p,q}$. By Theorem 2.6, we have $\per
L(C_4(1^02,1^{q-3}0,1^{p-3}0,1^00)) < \per L(G)$, a
contradiction to the choice of $G$. So $G\not\cong G_2$.

If $G\cong G_1$, by a simple calculation, we get
\begin{eqnarray}
 \per L(G_1) &=& 100(p-2)(q-2)+40n-140, \notag\\
  \per L(B(p,q)) &=& 20(p-1)(q-1)+4n.\label{eq:4.9}
\end{eqnarray}

Note that $p+q=n$ with $3\le p\le q\le n-3$, hence
\[\label{eq:3.17}
  pq\geq 3(n-3).
\]
Hence,
\begin{eqnarray*}
 \per L(G_1)-\per L(B(p,q)) &=& 80pq-144n+240 \\
   &\geq &  80\cdot 3 (n-3)-144n+240  \ \ \ \ \text{(by (\ref{eq:3.17}))}\\
    &=& 96n-480 \\
    &>& 0.
\end{eqnarray*}
Therefore $\per L(G_1)> \per L(B(p,q))$. Thus we obtain $G \cong
B(p,q)$. Together with Eq. (\ref{eq:4.9}),  we complete the proof.
\end{proof}
\begin{figure}[h!]
\begin{center}
\psfrag{a}{$n-7$} \psfrag{b}{$n-6$} \psfrag{c}{$n-5$}
\psfrag{1}{$\hat{G}_1$} \psfrag{2}{$\hat{G}_2$} \psfrag{3}{$\hat{G}_3$}
\psfrag{4}{$\hat{G}_4$} \psfrag{5}{$\hat{G}_5$} \psfrag{6}{$\hat{G}_6$}
\psfrag{7}{$\hat{G}_7$} \psfrag{8}{$\hat{G}_8$}
\includegraphics[width=110mm]{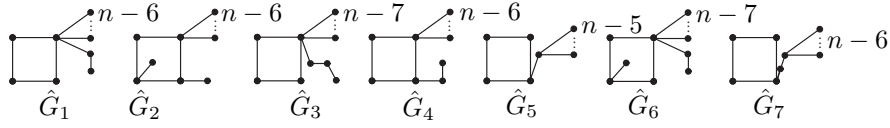}\\
\caption{Graphs $\hat{G}_1, \hat{G}_2, \ldots, \hat{G}_7$.}
\end{center}
\end{figure}

Next we are to identify the graph in $\mathscr{U}_{n}^{p,q}$ with the second (resp. third) smallest Laplacian permanent.
\begin{thm}
Among $\mathscr{U}_{n}^{p,q}.$
\begin{wst}
  \item[{\rm (i)}] If $p =2$, the ordering of all the members in $\mathscr{U}_n^{2,n-2}$ with $n\geq 4$ is as follows:
$
\per L(B(2,n-2))< \per L(C_4(1^10,1^00,1^{n-5}0,1^00))
< \cdots < \per L(C_4(1^i0,1^00, 1^{n-4-i}0,1^00))<\cdots< \per L(C_4(1^{\lfloor \frac{n-4}{2}\rfloor}0,1^00,\linebreak 1^{\lceil \frac{ n-4 }{2}\rceil}0,1^00)).
$\vspace{3mm}
\item[{\rm (ii)}] If $p=3,$ then $\per L(B(3,n-3))<\per L(\hat{G}_1)<\per L(\hat{G}_2)<\per L(G)$ for all $G\in \mathscr{U}_n^{3,n-3}\setminus\{B(3,n-3), \hat{G}_1, \hat{G}_2\}$ with $n\ge 20,$ where $\hat{G}_1, \hat{G}_2$ are depicted in Fig. 8.
  \item[{\rm (iii)}] If $p\geq 4,$
  \begin{wst}
    \item[{\rm (a)}] for all $G\in \mathscr{U}_{n}^{p,q} \setminus \{B(p,q)\}$ with $n \ge 8$, one has
$ \per L(G)\geq 36(p-2)(q-1)+4p+8q-4$ with equality if and
only if $G\cong C_4(1^{q-2}0, 1^{p-3}0,1^00,1^10)$.\vspace{3mm}
    \item[{\rm (b)}] for all $G \in \mathscr{U}_{n}^{p,q}
  \setminus \{B(p,q),C_4(1^{q-2}0,1^{p-3}0,1^00,1^10)\}$ with $q >p,\, n\geq 9$, one has
$ \per L(G)\geq 36(q-2)(p-1)+4q+8p-4$ with equality if and
only if $G\cong C_4(1^{q-3}0,1^{p-2}0,1^10,1^00).$
  \end{wst}
\end{wst}
\end{thm}
\begin{proof}
(i)\ If $p=2$, then
$$
\mathscr{U}_n^{2,n-2}=\left\{C_4(1^i2^0,1^02^0,1^{n-4-i}2^0,1^02^0):0\leq
i\leq \left\lfloor\frac{ n-4 }{2}\right\rfloor\right\}.
$$
By a simple calculation, we get
$$
\per L(C_4(1^i2^0,1^02^0,1^{n-4-i}2^0,1^02^0))
=-16\left(i-\frac{n-4}{2}\right)^2+4(n-1)^2.
$$
Consider the function $f(x)=-16(x- \frac{ n-4 }{2})^2+4(n-1)^2$ in $x$ with
$0\leq x \leq \lfloor\frac{ n-4 }{2}\rfloor.$ By the monotonicity of
$f(x)$, we have
\[\label{eq:4.8}
f(0)< f(1)< \cdots < f(i)<\cdots<f(\left\lfloor\frac{ n-4
}{2}\right\rfloor).
\]
Note that $C_4(1^00,1^00,1^{n-4}0,1^00)\cong B(2,n-2)$, hence (i)
follows immediately from (\ref{eq:4.8}).

(ii)\ Note that $p=3$, hence the cycle $C_r$ contained in $T\in
\mathscr{U}_n^{3,n-3}$ is of length at most 6, i.e., $r\le 6$. If
$r=4,$ the bipartite unicyclic graph, say $U$ (resp. $U'$), in
$\mathscr{U}_n^{3,n-3}$ with the second (resp. third) smallest
Laplacian permanent should satisfy the following property: Apply
Operation I (or I\!I) to $U$ (resp. $U'$) only once to get the graph
$B(3,n-3)$ (resp.  $B(3,n-3)$ or $U$). Hence, $U, U'\in
\{\hat{G}_{1}, \hat{G}_{2}, \hat{G}_{3}, \hat{G}_{4}, \hat{G}_5,
\hat{G}_6, \hat{G}_7, \hat{G}_8\}$, where $\hat{G}_1, \hat{G}_2,
\ldots, \hat{G}_8$ are depicted in Fig. 8.

If $r=6,$ then graph $G_1$ (see Fig. 7) is the possible graph with the smallest Laplacian permanent.
By direct calculation, we have
$$
\begin{array}{lll}
  \per L(\hat{G}_1) = 72n-276, & \per L(\hat{G}_2) = 76n-352,&\per L(\hat{G}_3) = 168n-804, \\
  \per L(\hat{G}_4) = 112n-516, &\per L(\hat{G}_5) = 96n-420, &\per L(\hat{G}_6) = 120n-580,\\
  \per L(\hat{G}_7) = 216n-1140,, &\per L(G_1) = 140n-640.  &
\end{array}
$$
Based on the above direct computing, (ii) follows immediately.

(iii)\ We first determine the graph, say $G$, in
$\mathscr{U}_n^{p,q}$ with the second smallest Laplacian permanent
for $p\ge 4$. In view of the proof of Theorem 3.7, it is easy to see
that the cycle $C_r$ contained in $G$ is of length at most $6.$
Furthermore, if $r=6$, only $G_1$ as depicted in Fig. 7 is
possible to be the particular graph $G$. If $r=4$, in view of
Theorems 2.2 and 2.4, we know that $B(p,q)$ can be obtained from $G$
by Operation I (or, I\!I) once. Hence, based on Operation I, $G$ may
be $C_4(1^{q-2}0,\,1^{p-3}0,1^00,1^10)$ or
$C_4(1^{q-3}0,1^{p-2}0,1^10,1^00);$ whereas based on Operation I\!I,
$G$ may be in the set
$$
\mathscr{A}=\{C_4(1^{q-3}(t+1), 1^{p-2-t}0, 1^00, 1^00): 1\le t\le p-2\}
$$
or
$$
\mathscr{B}=\{C_4(1^{q-2-t}0, 1^{p-3}(t+1), 1^00, 1^00): 1\le t\le q-2\}.
$$
Combining with Operation I we see that $\mathscr{A}$ (resp.
$\mathscr{B}$) contains just two members $C_4(1^{q-3}2, 1^{p-3}0,1^00,1^00)$ and $C_4(1^{q-3}(p-1), 1^00, 1^00,1^00)$ (resp. $C_4(1^{q-3}0, 1^{p-3}2, 1^00,1^00)$ and $C_4(1^00, 1^{p-3}(q-1),1^00,1^00)$). Hence,
summarizing the discussion as above we get that $G$ must be in $
\mathscr{U'}=\{C_4(1^{q-2}0, 1^{p-3}0,1^00,1^10),
C_4(1^{q-3}0,1^{p-2}0, 1^10,1^00),\linebreak C_4(1^{q-3}2, 1^{p-3}0,
1^00, 1^00), C_4(1^{q-3}0, 1^{p-3}2, 1^00,
1^00),$ $C_4(1^00, 1^{p-3}(q-1), 1^00, 1^00),
C_4(1^{q-3}(p-1), 1^00, 1^00, 1^00),\linebreak G_1\}.$

By direct calculation, we obtain
\begin{eqnarray}
  \per L(C_4(1^{q-2}0, 1^{p-3}0,1^00,1^10)) &=& 36pq-32n-32q+68,\notag \\
  \per L(C_4(1^{q-3}0, 1^{p-2}0,1^10,1^00)) &=& 36pq-32n-32p+68,\notag \\
 \per L(C_4(1^{q-3}2, 1^{p-3}0,1^00,1^00)) &=& 60pq-68n-40q+144,\notag \\
 \per L(C_4(1^{q-3}0, 1^{p-3}2,1^00,1^00)) &=& 60pq-68n-40p+144,\notag \\
 \per L(C_4(1^{q-3}(p-1), 1^00,1^00,1^00)) &=& 48pq-72n+24p+84,\notag \\
 \per L(C_4(1^00, 1^{p-3}(q-1),1^00,1^00)) &=& 48pq-72n+24q+84.\notag
\end{eqnarray}
This gives
\[\label{eq:3.9}
\per L(C_4(1^{q-2}0,1^{p-3}0,1^00,1^10)) < \per
L(C_4(1^{q-3}0,1^{p-2}0,1^10,1^00))< \per L(\hat{G})
\]
for all $\hat{G}\in \mathscr{U}' \setminus \{
C_4(1^{q-2}0,1^{p-3}0,1^00,1^10),
C_4(1^{q-3}0,1^{p-2}0,1^10,1^00)\}$ for $q > p \geq 4$. This
completes the proof of the first part of (iii).

Now we show the second part of (iii). By a similar discussion as in
the proof of the first part of (iii), we know that the graph, say
$G'$, in $\mathscr{U}_{n}^{p,q}$ having the third smallest Laplacian
permanent is either the graph with the second smallest Laplacian
permanent in $\mathscr{U}'$, or apply Operation I (or I\!I) once to
$G'$ to obtain the graph $C_4(1^{q-2}0,1^{p-3}0,1^00,1^10)$, which
has the second smallest Laplacian permanent in
$\mathscr{U}_{n}^{p,q}$. Hence, together with (\ref{eq:3.9}), we
obtain that $G'$ is in the set $ \mathscr{U}''
=\{C_4(1^{q-3}0, 1^{p-2}0, 1^10,1^00),
C_4(1^{q-3}2,1^{p-3}0,1^00,1^00),C_4(1^{q-3}0,
\linebreak 1^{p-4}2, 1^00,1^10),C_4(1^{q-3}2,
1^{p-4}0, 1^00,\,1^10),C_4(1^{q-3}0,\, 1^{p-3}0,\, 1^00,\, 1^02),
C_4(1^{q-3}0, 1^{p-3}0,\,1^10,\,1^10), $ $C_4(1^00,
1^{p-4}(q-1),1^00,1^10),\,C_4(1^{q-3}(p-2),1^00,1^00,1^10),\,C_4(1^00, 1^{p-3}0,1^00,1^0(q-1))\}.$

By direct calculation, we have
\begin{eqnarray}
  \per L(C_4(1^{q-3}2,1^{p-3}0,1^00,1^00)) &=& 60pq-68n-40q+144, \notag\\
  \per L(C_4(1^{q-3}0, 1^{p-4}2,1^00,1^10)) &=& 108pq-24q -204n+464,\notag \\
  \per L(C_4(1^{q-3}2, 1^{p-4}0,1^00,1^10)) &=& 108pq-168q-132n+400,\notag \\
 \per L(C_4(1^{q-3}0, 1^{p-3}0,1^00,1^02)) &=& 92pq+8q-172n+336,\notag \\
\per L(C_4(1^{q-3}0, 1^{p-3}0,1^10,1^10)) &=& 68pq-128n+260,\notag \\
 \per L(C_4(1^00, 1^{p-4}(q-1),1^00,1^10)) &=& 80pq-40q-120n+260,\notag \\
 \per L(C_4(1^{q-3}(p-2), 1^00,1^00,1^10)) &=& 88pq-120q-100n+272,\notag \\
 \per L(C_4(1^00, 1^{p-3}0,1^00,1^0(q-1))) &=& 64pq+8p-96n+132.\notag
\end{eqnarray}
Based on the above direct computing, the second part of (iii) follows immediately.
\end{proof}

\noindent{\bf Remark\ 3.} In view of Theorems 3.7 and 3.8, we hope
to show that, among the set of all $n$-vertex unicyclic graphs with
a $(p,q)$-bipartition($q> p\geq 4$),  $B(p,q) \prec
C_4(1^{q-2}0,1^{p-3}0, 1^00,1^10) \prec C_4(1^{q-3}0, 1^{p-2}0, 1^10,1^00) \prec G$
for all $G\in \mathscr{U}_n^{p,q}\setminus  \{B(p,q), C_4(1^{q-2}0,
1^{p-3}0, 1^00,1^10), C_4(1^{q-3}0, 1^{p-2}0, 1^10,1^00)\}$ in the
future research. If this is true, it will support the relationship
between the Laplacian coefficients and the Laplacian permanent of
$n$-vertex bipartite unicyclic graphs with a $(p,q)$-bipartition.

To conclude this subsection, we determine the first, second, third
smallest Laplacian permanent of graphs in $\mathscr{U}_n$, the set
of all bipartite unicyclic graphs on $n$ vertices.
\begin{thm}
Among $\mathscr{U}_n$ with $n\geq 4$,
\begin{wst}
  \item[{\rm (i)}] for all $G\in \mathscr{U}_n,$ we have
$
\per L(G)\geq 24n-60
$
with equality if and only if $G\cong B(2,n-2)$.
  \item[{\rm (ii)}] for all $G \in \mathscr{U}_n
  \setminus \{B(2,n-2)\}$ with $n\geq 6$, we have
$
\per L(G)\geq 40n-140
$
with equality if and only if $G\cong C_4(1^10,1^00,1^{n-5}0,1^00).$
\item[{\rm (iii)}] for all $G \in \mathscr{U}_n \setminus \{B(2,n-2),C_4(1^10,1^00,1^{n-5}0,1^00)\}$ with $n\geq 6$, we have
$
\per L(G)\geq 44n-160
$
with equality if and only if $G\cong B(3,n-3).$
\end{wst}
\end{thm}
\begin{proof}
It is routine to see that $\mathscr{U}_n=\mathscr{U}_n^{2,n-2} \cup
\mathscr{U}_n^{3,n-3} \cup \cdots \cup \mathscr{U}_n^{\lfloor
\frac{n}{2} \rfloor, \lceil \frac{n}{2}\rceil}$. Note that for all
$G\in \mathscr{U}_n^{p,q}$, by Theorem 3.7 one has $\per L(G)\geq
\per L(B(p,q))=20(p-1)(q-1)+4n,$ with the equality if and only if
$G\cong B(p,q).$ Consider the function
$$
f(x)=20(x-1)(n-x-1)+4n
$$
in $x$ with $2\le x\le \lfloor\frac{n}{2}\rfloor.$ It is routine to
check that $f'(x)=20(n-2x)>20(n-x-(n-x))=0$. Hence, $f(x)$ is an
increasing function for $2\le x\le \lfloor\frac{n}{2}\rfloor.$ That
is to say, $f(2)<f(3)<\cdots<f(\lfloor\frac{n}{2}\rfloor),$ which
implies (i) immediately.

Based on Theorems 3.7-3.8 and the proof in (i) as above, in order to
determine the the graph in $\mathscr{U}_n$ having the second minimal
Laplacian permanent, it suffices to compare the values between $\per
L(C_4(1^10,1^00,1^{n-5}0,1^00))$ and $\per L(B(3,n-3)).$ By an
elementary calculation, we have
\[
  \per L(C_4(1^10,1^00,1^{n-5}0,1^00)) = 40n-140,\ \ \  \per L(B(3,n-3))=44n-160.
\]
It is routine to check that $\per L(C_4(1^10,1^00,1^{n-5}0,1^00)) <
\per L(B(3,n-3)).$ Hence, (ii) holds immediately.

Similarly, in order to determine the third minimal Laplacian
permanent among $\mathscr{U}_n$, it suffices to compare the values
between $\per L(C_4(1^20,1^00,1^{n-6}0,1^00))$ and $\per
L(B(3, n-3))$. Note that if $n=6$ (resp. $7$), it is
straightforward to check that
$C_4(1^20,1^00,1^{n-6}0,1^00)$ does not exist and
$B(3,n-3)$ is the graph with the third minimal Laplacian permanent
among $\mathscr{U}_n$. For $n\ge 8$, by direct calculation, we have
\[
\per L(C_4(1^20,1^00,1^{n-6}0,1^00)) = 56n-252.
\]
In view of the second equation in (3.10) and (3.11), it is routine
to check that $ \per L(C_4(1^20,1^00,1^{n-6}0,1^00))>\per
L(B(3,n-3))=44n-160.$ Hence, (iii) holds immediately.

This completes the proof.
\end{proof}
\section*{\normalsize Acknowledgements}
The authors would like to express their sincere gratitude to the
referee for a very careful reading of the paper and for all his or
her insightful comments and valuable suggestions, which led to a
number of improvements in this paper.

\end{document}